\documentclass[11pt]{amsart}
\usepackage[utf8]{inputenc}
\usepackage{amsfonts,amsmath,amssymb,amsthm,hyperref,setspace}
\usepackage[a4paper,margin=1.9cm,top=2.5cm,bottom=2.5cm,centering,vcentering]{geometry}

\usepackage[normalem]{ulem}

\newcommand{\op}{\overline{p}}
\usepackage{color}

\newtheorem{theorem}{Theorem}
\newtheorem{conjecture}{Conjecture}

\newtheorem{remark}[theorem]{Remark}

\title{Some Missed Congruences modulo powers of $2$ for $t$-colored overpartitions}

\author{Manjil P. Saikia}
\address{School of Mathematics, Cardiff University, Cardiff CF24 4AG, United Kingdom}
\address{Current Address: Department of Humanities and Basic Sciences, Indian Institute of Information Technology (IIIT) Manipur, Imphal 795002, Manipur, India}
\email{manjil@saikia.in}

\keywords{integer partitions, Ramanujan-type congruences, Radu's algorithm.fb}

\subjclass[2020]{11P81, 11P83.}

\begin{document}

\maketitle

\begin{abstract}
    Recently, Nayaka and Naika (2022) proved several congruences modulo $16$ and $32$ for $t$-colored overpartitions with $t=5,7,11$ and $13$. We extend their list using an algorithmic technique.
\end{abstract}

\section{Introduction}

A partition of a positive integer $n$ is a finite non-increasing sequence of positive integers $\lambda=(\lambda_1, \lambda_2, \ldots, \lambda_k)$ such that $\sum\limits_{i=1}^k\lambda_i=n$. The number of partitions of $n$ is denoted by $p(n)$. Euler found the generating function of $p(n)$, given by
\[
\sum_{n\geq 0}p(n)q^n=\frac{1}{(q;q)_\infty},
\]
where
\[
(a;q)_\infty:=\prod_{i\geq 0}(1-aq^i), \quad |q|<1.
\]

Ramanujan found several interesting congruences modulo primes that the partition function satisfies, and this opened up the area to further exploration. Several mathematicians have studied the arithmetic properties of the partition function as well as other generalized classes of partitions. In this tradition, recently Nayaka and Naika \cite{NayakaNaika} studied the congruences modulo $2$ safisfied by the $t$-colored overpartitions.

An overpartition of a nonnegative integer $n$ is a non-increasing sequence of natural numbers whose sum is $n$, and where the first occurrence (equivalently, the last occurrence) of a number may be overlined. The number of overpartitions of $n$ is denoted by $\op(n)$ and its generating function is given by
\[
\sum_{n\geq 0}\op(n)q^n=\frac{(q^2;q^2)_\infty}{(q;q)^2_\infty}.
\]
Further, a partition is called a $t$-colored partition if each part can appear as $t$ colors. The number of $t$-color partitions of $n$ is denoted by $p_{-t}(n)$ and its generating function is given by
\[
\sum_{n\geq 0}p_{-t}(n)q^n=\frac{1}{(q;q)^t_\infty}.
\]
Now, we define the function $\op_{-t}(n)$ whose generating function is given by
\begin{equation}\label{eq:eob}
    \sum_{n\geq 0}\op_{-t}(n)q^n=\frac{(q^2;q^2)^t_\infty}{(q;q)^{2t}_\infty}.
\end{equation}
This function denotes the number of $t$-colored overpartitions of $n$.

Nayaka and Naika \cite{NayakaNaika} used elementary techniques to find several congruences modulo $16$ and $32$ satisfied by $\op_{-t}(n)$ for $t=5,7,11$ and $13$. The goal of the present short note is to prove algorithmically several missed congruences modulo powers of $2$ satisfied by $\op_{-t}(n)$ for $t=5,7,11$ and $13$, in certain cases bettering the congruences found by Nayaka and Naika \cite{NayakaNaika}. We state our results below.

\begin{theorem}\label{thm:5}
For all $n\geq 0$, we have
\begin{align}
  \op_{-5}(8n+1) &\equiv 0 \pmod{2},\\
 \op_{-5}(8n+2) &\equiv 0 \pmod{4},\\
 \op_{-5}(8n+3) &\equiv 0 \pmod{8},\\
 \op_{-5}(8n+4) &\equiv 0 \pmod{2},\\
 \op_{-5}(8n+5) &\equiv 0 \pmod{8},\\
  \op_{-5}(8n+6) &\equiv 0 \pmod{8},\\
 \op_{-5}(8n+7) & \equiv 0 \pmod{128}.
\end{align}
\end{theorem}

\begin{remark}
Nayaka and Naika \cite{NayakaNaika} had obtained
\[
\op_{-5}(8n+7)  \equiv 0 \pmod{32}.
\]
\end{remark}

\begin{theorem}\label{thm:7}
For all $n\geq 0$, we have
\begin{align}
  \op_{-7}(8n+1) &\equiv 0 \pmod{2},\\
 \op_{-7}(8n+2) &\equiv 0 \pmod{16},\\
 \op_{-7}(8n+3) &\equiv 0 \pmod{32},\\
 \op_{-7}(8n+4) &\equiv 0 \pmod{2},\\
 \op_{-7}(8n+7) & \equiv 0 \pmod{128}.
\end{align}
\end{theorem}

\begin{theorem}\label{thm:11}
For all $n\geq 0$, we have
\begin{align}
  \op_{-11}(8n+1) &\equiv 0 \pmod{2},\\
 \op_{-11}(8n+2) &\equiv 0 \pmod{8},\\
 \op_{-11}(8n+3) &\equiv 0 \pmod{16},\\
 \op_{-11}(8n+4) &\equiv 0 \pmod{2},\\
 \op_{-11}(8n+7) & \equiv 0 \pmod{64}.
\end{align}
\end{theorem}

\begin{remark}
Nayaka and Naika \cite{NayakaNaika} had obtained
\[
\op_{-11}(8n+7)  \equiv 0 \pmod{32}.
\]
\end{remark}

\begin{theorem}\label{thm:13}
For all $n\geq 0$, we have
\begin{align}
  \op_{-13}(8n+1) &\equiv 0 \pmod{2},\\
 \op_{-13}(8n+2) &\equiv 0 \pmod{4},\\
 \op_{-13}(8n+3) &\equiv 0 \pmod{8},\\
 \op_{-13}(8n+4) &\equiv 0 \pmod{2},\\
 \op_{-13}(8n+5) &\equiv 0 \pmod{8},\\
  \op_{-13}(8n+6) &\equiv 0 \pmod{8},\\
 \op_{-13}(8n+7) & \equiv 0 \pmod{256}.
\end{align}
\end{theorem}

\begin{remark}
Nayaka and Naika \cite{NayakaNaika} had obtained
\[
\op_{-13}(8n+7)  \equiv 0 \pmod{32}.
\]
\end{remark}

Theorems \ref{thm:5}, \ref{thm:7}, \ref{thm:11} and \ref{thm:13} can be proved using Smoot's \cite{Smoot} implementation of an algorithm of Radu \cite{Radu} which we will describe in the next section. The above results suggest that several more congruences might be true. We make the following conjecture.

\begin{conjecture}
For all $n\geq 0$ and primes $q$, we have
\begin{align}
\op_{-q}(8n+1) & \equiv 0 \pmod{2},\\
\op_{-q}(8n+2) & \equiv 0 \pmod{4},\\
\op_{-q}(8n+3) & \equiv 0 \pmod{8},\\
\op_{-q}(8n+4) & \equiv 0 \pmod{2},\\
\op_{-q}(8n+5) & \equiv 0 \pmod{8},\\
\op_{-q}(8n+6) & \equiv 0 \pmod{8},\\
\op_{-q}(8n+7) & \equiv 0 \pmod{32}.
\end{align}
\end{conjecture}

Nayanka and Naika \cite{NayakaNaika} also proved several infinite families of congruences modulo $16$ and $32$ for $\op_{-t}(n)$ with $t=7,11$ and $13$. Using our methods combined with elementary reasoning we can also prove such infinite families. As an example we give the following family which was not proved by Nayaka and Naika.

\begin{theorem}\label{thm:inf}
For all $\alpha, \beta, \gamma \geq 0$, we have
\begin{equation}\label{eq:inf}
    \sum_{n\geq 0}\op_{-5}(8\cdot 3^{2\alpha}\cdot 5^{2\beta}\cdot 7^{2\gamma}n+2\cdot 3^{2\alpha}\cdot 5^{2\beta}\cdot 7^{2\gamma})\equiv 4f_1^6 \pmod{8},
\end{equation}
\begin{equation}\label{eq:inf2}
    \sum_{n\geq 0}\op_{-5}(8\cdot 3^{2\alpha+1}\cdot 5^{2\beta}\cdot 7^{2\gamma}n+2\cdot 3^{2\alpha+2}\cdot 5^{2\beta}\cdot 7^{2\gamma})\equiv 4f_3^6 \pmod{8},
\end{equation}
\begin{equation}\label{eq:inf3}
    \sum_{n\geq 0}\op_{-5}(8\cdot 3^{2\alpha}\cdot 5^{2\beta+1}\cdot 7^{2\gamma}n+2\cdot 3^{2\alpha}\cdot 5^{2\beta+1}\cdot 7^{2\gamma})\equiv 4qf_5^6 \pmod{8},
\end{equation}
and
\begin{equation}\label{eq:inf4}
    \sum_{n\geq 0}\op_{-5}(8\cdot 3^{2\alpha}\cdot 5^{2\beta}\cdot 7^{2\gamma+1}n+2\cdot 3^{2\alpha}\cdot 5^{2\beta}\cdot 7^{2\gamma+1})\equiv 4qf_7^6 \pmod{8},
\end{equation}
where $f_k:=(q^k;q^k)_\infty$.
\end{theorem}
\noindent We prove this result in Section \ref{sec:inf} using a combination of algorithmic and elementary techniques.

\section{Proofs of Theorems \ref{thm:5}, \ref{thm:7}, \ref{thm:11} and \ref{thm:13} }

Theorems \ref{thm:5}, \ref{thm:7}, \ref{thm:11} and \ref{thm:13}  can be proved using Smoot's \cite{Smoot} implementation of an algorithm of Radu \cite{Radu}. Radu's algorithm can be used to prove Ramanujan type congruences of the form stated in the previous section. The algorithm takes as an input the generating function
\[
\sum_{n\geq 0}a_r(n)q^n=\prod_{\delta|M}\prod_{n\geq 1}(1-q^{\delta n})^{r_\delta},
\]
and positive integers $m$ and $N$, with $M$ another positive integer and $(r_\delta)_{\delta|M}$ is a sequence indexed by the positive divisors $\delta$ of $M$. With this input, Radu's algorithm tries to produce a set $P_{m,j}(j)\subseteq \{0,1,\ldots, m-1\}$ which contains $j$ and is uniquely defined by $m, (r_\delta)_{\delta|M}$ and $j$. Then, it decides if there exists a sequence $(s_\delta)_{\delta |N}$ such that
\[
q^\alpha \prod_{\delta|M}\prod_{n\geq 1}(1-q^{\delta n})^{s_\delta} \cdot \prod_{j^\prime \in P_{m,j}(j)}\sum_{n\geq 0}a(mn+j^\prime)q^n,
\]
is a modular function with certain restrictions on its behaviour on the boundary of the upper half complex plane $\mathbb{H}$.

Smoot \cite{Smoot} implemented this algorithm in Mathematica and we will use his \texttt{RaduRK} package which requires the software packaage \texttt{4ti2}. Documentation on how to install and use these packages are available from Smoot \cite{Smoot}. Since the proofs of all the results are similar, we are only going to prove one equation from Theorem \ref{thm:5}, namely
\begin{equation}\label{eq:pp}
    \op_{-5}(8n+7)\equiv 0 \pmod{128}.
\end{equation}

It is natural to guess that $N=m$ (which corresponds to the congruence subgroup $\Gamma_0(N)$), but this is not always the case, although they are usually closely related to one another. The determination of the correct value of $N$ is an important problem for the usage of \texttt{RaduRK} and it depends on a criterion called the $\Delta^\ast$ criterion \cite[Definitions 34 and 35]{Radu}, which we do not explain here. It is easy to check the minimum $N$ which satisfies this criterion by running \texttt{minN[M, r, m, j]} which we do now for our generating functions. The generating function of $\op_{-t}(n)$ given in \eqref{eq:eob} can be described by setting $M=2$ and $r=\{-10,5\}$. Now, running \texttt{minN[2, \{-10,5\}, 8, 7]} yields $8$ as the output for the minimum choice of $N$.

Running \texttt{RK[8,2,{-10,5},8,7]} gives the following.
  \allowdisplaybreaks{
	\begin{doublespace}
		\begin{align*}
		\texttt{In[1] := } & \texttt{RK[8,2,\{-10,5\},8,7]}\\
		& \prod_{\texttt{$\delta$|M}} (\texttt{q}^{\delta };\texttt{q}^{\delta })_{\infty }^{\texttt{r}_{\delta }}  = \sum_{\texttt{n=0}}^{\infty } \texttt{a}(\texttt{n})\,\texttt{q}^\texttt{n}\\
		& \fbox{$\texttt{f}_\texttt{1}(\texttt{q})\cdot \prod\limits_{\texttt{j}'\in \texttt{P}_{\texttt{m,r}}(\texttt{j}) } \sum\limits_{\texttt{n=0}}^\infty \texttt{a}(\texttt{mn}+\texttt{j}')\,\texttt{q}^\texttt{n} = \sum\limits_{\texttt{g}\in \texttt{AB}} \texttt{g}\cdot \texttt{p}_\texttt{g}(\texttt{t}) $} \\
		& \texttt{Modular Curve: }\texttt{X}_\texttt{0}(\texttt{N}) \\
		\texttt{Out[2] = }\\
& \begin{array}{c|c}
 \text{N:} & 8 \\
\hline
 \text{$\{$M,(}r_{\delta })_{\delta |M}\text{$\}$:} & \{2,\{-10,5\}\} \\
\hline
 \text{m:} & 8 \\
\hline
 P_{m,r}\text{(j):} & \{7\} \\
\hline
 f_1\text{(q):} & \dfrac{(q;q)_{\infty }^{79} \left(q^4;q^4\right)_{\infty }^{36}}{q^{17}
   \left(q^2;q^2\right)_{\infty }^{38} \left(q^8;q^8\right)_{\infty }^{72}} \\
\hline
 \text{t:} & \dfrac{\left(q^4;q^4\right)_{\infty }^{12}}{q \left(q^2;q^2\right)_{\infty
   }^4 \left(q^8;q^8\right)_{\infty }^8} \\
\hline
 \text{AB:} & \{1\} \\
\hline
 \left\{p_g\text{(t): g$\in $AB$\}$}\right. & \text{given below in \eqref{eq:pgt}} \\
\hline
 \text{Common Factor:} & 128 \\
\end{array}
		\end{align*}
	\end{doublespace}  }

\noindent where
\begin{multline}\label{eq:pgt}
   \left\{p_g\text{(t): g$\in $AB$\}$}\right.=  \left\{37760 t^{17}+47761408
   t^{16}+10846240768 t^{15}\right. +868870094848 t^{14}\\+32519056130048 t^{13}+661947909931008
   t^{12}+8026570602053632 t^{11}\\+61243801104023552 t^{10}+302871878945472512
   t^9+979900817664376832 t^8\\+2054802074125729792 t^7+2711338639077408768
   t^6+2131168862538825728 t^5\\ \left.+911076328575336448 t^4+181969724152741888
   t^3+12820855335682048 t^2+162177965096960 t\right\}.
\end{multline}
\noindent This shows that
\[
f_1(q)\cdot \Bigg(\sum_{n\geq 0}\op_{-5}(8n+7)q^n\Bigg)=p_g(t),
\]
where $f_1(q)=\dfrac{(q;q)_{\infty }^{79} \left(q^4;q^4\right)_{\infty }^{36}}{q^{17}
   \left(q^2;q^2\right)_{\infty }^{38} \left(q^8;q^8\right)_{\infty }^{72}}$ and $t=\dfrac{\left(q^4;q^4\right)_{\infty }^{12}}{q \left(q^2;q^2\right)_{\infty
   }^4 \left(q^8;q^8\right)_{\infty }^8}$. This immediately shows that equation \eqref{eq:pp} is true, as the last row of the output gives the common factor.

   The rest of the equations in Theorems \ref{thm:5}, \ref{thm:7}, \ref{thm:11} and \ref{thm:13} can be proved in exactly a similar way, so we omit the details here. The interested reader can download the Mathematica file from \url{https://manjilsaikia.in/publ/opt.nb} to check the relevant outputs.

\section{Proof of Theorem \ref{thm:inf}}\label{sec:inf}

It can be proved using \texttt{RaduRK} that
\begin{equation}\label{eq:1}
    \sum_{n\geq 0}p_{-5}(8n+2)q^n\equiv 4\frac{f_4^{179}}{f_1^{78}f_2^{36}f_8^{70}}\pmod{8}.
\end{equation}
We do not show the details here. The interested reader can download the Mathematica file from \url{https://manjilsaikia.in/publ/opt-gf.nb} to check the relevant output. Using the congruence
\[
f_{m}^{2^k}\equiv f_{2m}^{2^{k-1}}\pmod {2^k},
\]
we arrive from equation \eqref{eq:1} at
\begin{equation}
   \sum_{n\geq 0}p_{-5}(8n+2)q^n\equiv  4f_1^6 \pmod{8}.
\end{equation}
This is the case $\alpha=\beta=\gamma=0$ of \eqref{eq:inf}.

Let us now consider the case $\beta=\gamma=0$. We need the following formula \cite[p. 345, Entry 1(iv)]{BerndtIII}
\begin{equation}\label{11}
    f_1^3=\frac{f_6f_9^3}{f_3f_{18}^3}+4q^3\frac{_3^2f_{18}^6}{_6^2f_9^3}-3qf_9^3\equiv f_3+qf_9^3 \pmod 2.
\end{equation}
Using this in equation \eqref{eq:inf} with $\beta=\gamma=0$ we have
\begin{equation*}
    \sum_{n\geq 0}\op_{-5}(8\cdot 3^{2\alpha}n+2\cdot 3^{2\alpha})q^n\equiv 4q^2f_9^6+4f_3^2 \pmod 8,
\end{equation*}
for all $n\geq 0$. Extracting the terms involving $q^{3n+2}$ from the above gives us
\begin{equation*}
    \sum_{n\geq 0}\op_{-5}(8\cdot 3^{2\alpha+1}n+2\cdot 3^{2\alpha+2})q^n\equiv 4f_3^6 \pmod 8,
\end{equation*}
which gives us
\begin{equation*}
    \sum_{n\geq 0}\op_{-5}(8\cdot 3^{2\alpha+2}n+2\cdot 3^{2\alpha+2})q^n\equiv 4f_1^6 \pmod 8.
\end{equation*}
This implies that equation \eqref{eq:inf} is true for $\alpha+1$. Hence by induction it is true for any non-negative integer $\alpha$ and $\beta=\gamma=0$.

Now we consider the case $\gamma=0$ and let equation \eqref{eq:inf} be true for some integer $\alpha, \beta\geq 0$. We have the following formula \cite[Theorem 7.4.4]{Spirit}
\begin{align}\label{13}
f_1 = f_{25} \Bigg(\frac{1}{R(q^5)} - q - q^2 R(q^5)\Bigg),
\end{align}
where
\[
R(q)=\frac{(q;q^5)_\infty (q^4;q^5)_\infty}{(q^2;q^5)_\infty(q^3;q^5)_\infty}.
\]
Using the above in equation \eqref{eq:inf} we ontain
\begin{equation*}
    \sum_{n\geq 0}\op_{-5}(8\cdot 3^{2\alpha}\cdot 5^{2\beta}n+2\cdot 3^{2\alpha}\cdot 5^{2\beta})q^n\equiv 4f_{25}^6\Bigg(\frac{1}{R(q^5)} - q - q^2 R(q^5)\Bigg)^6 \pmod 8.
\end{equation*}
Now we extract the terms involving $q^{5n+1}$ from the above to arrive at
\begin{equation*}
        \sum_{n\geq 0}\op_{-5}(8\cdot 3^{2\alpha}\cdot 5^{2\beta+1}n+2\cdot 3^{2\alpha}\cdot 5^{2\beta+2})q^n\equiv 4qf_{5}^6 \pmod 8,
\end{equation*}
which implies
\begin{equation*}
        \sum_{n\geq 0}\op_{-5}(8\cdot 3^{2\alpha}\cdot 5^{2\beta+2}n+2\cdot 3^{2\alpha}\cdot 5^{2\beta+2})q^n\equiv 4f_{1}^6 \pmod 8,
\end{equation*}
which implies equation \eqref{eq:inf} is true for $\beta+1$ and hence  by indeuction it is true for all $\alpha, \beta\geq 0$.

Now we suppose that the equation \eqref{eq:inf} holds for some integers $\alpha, \beta, \gamma \geq 0$. We need the following $7$-dissection formula \cite[p. 303, Entry 17(v)]{BerndtIII} for $(q;q)_\infty$
\begin{equation}\label{14}
 f_1=f_{49} (A^\prime(q^7)-qB^\prime(q^7)-q^2+q^5C^\prime(q^7)),
\end{equation}
where
\[
A^\prime(q^7)=\frac{f(-q^{14},-q^{35})}{f(-q^7,-q^{42})}, \quad B^\prime(q^7)=\frac{f(-q^{21},-q^{28})}{f(-q^{14},-q^{35})}, \quad \text{and} \quad C^\prime(q^7)=\frac{f(-q^7,-q^{42})}{f(-q^{21},-q^{28})},
\]
and
\[
f(a,b)=\sum_{n=-\infty}^\infty a^{n(n+1)/2}b^{n(n-1)/2}, \quad |ab|<1.
\]
In a similar fashion like earlier, using the above dissection formula and extracting the terms involving $q^{7n+5}$ and then extracting the coefficient of $q^{7n+1}$ from the resulting congruence we shall obtain
\[
\sum_{n\geq 0}\op_{-5}(8\cdot 3^{2\alpha}\cdot 5^{2\beta}\cdot 7^{2\gamma+2}n+2\cdot 3^{2\alpha}\cdot 5^{2\beta}\cdot 7^{2\gamma+2})q^n \equiv 4f_1^6 \pmod 8,
\]
which implies equation \eqref{eq:inf} is true for $\gamma+1$ and hence it is true for any non-negative integers $\alpha, \beta, \gamma\geq 0$. Thus, we have proved equation \eqref{eq:inf}.

We can prove equation \eqref{eq:inf2} using equation \eqref{11} in equation \eqref{eq:inf}; equation \eqref{eq:inf3} using equation \eqref{13} in equation \eqref{eq:inf} and equation \eqref{eq:inf4} using equation \eqref{14} in equation \eqref{eq:inf}.

\section{Concluding Remarks}\label{sec:conc}

Equation \eqref{14} is a special case of a result of Ramanathan \cite[Theorem 1]{Ramanathan} (also independently proved by Evans \cite{Evans}) which is stated below.

\begin{theorem}\cite[Theorem 12.1]{BerndtIII}
    Let $n$ be a natural number with $n\equiv 1 \pmod 6$ and let $g\geq 1$. If $n=6g+1$, then
    \begin{equation}
    f_1= f_{n^2}\Bigg((-1)^gq^{(n^2-1)/24}+\sum_{k=1}^{(n-1)/2}(-1)^{k+g}q^{(k-g)(3k-3g-1)/2}\frac{f(-q^{2nk},-q^{n^2-2nk})}{f(-q^{nk},-q^{n^2-nk})} \Bigg).
    \end{equation}
    If $n=6g-1$, then
    \begin{equation}
  f_1= f_{n^2}\Bigg((-1)^gq^{(n^2-1)/24}+\sum_{k=1}^{(n-1)/2}(-1)^{k+g}q^{(k-g)(3k-3g+1)/2}\frac{f(-q^{2nk},-q^{n^2-2nk})}{f(-q^{nk},-q^{n^2-nk})} \Bigg). 
    \end{equation}
\end{theorem}
\noindent It might be possible to prove more congruences using this result like we used equation \eqref{14} in the previous section.

As can be guessed, Radu's algorithm will yield several more congruences. It would be interesting to see if there exists any general pattern. In particular we can prove a result of the type of Theorem \ref{thm:inf} for $\op_{-13}(n)$ as well which was missed by Nayaka and Naika \cite{NayakaNaika}. We leave the details to the reader.

\subsection*{Funding}

The author is partially supported by the Leverhulme Trust Research Project Grant RPG-2019-083.

\end{document}